\documentclass[a4paper,12pt]{article}

\usepackage{amsthm}
\usepackage{amsmath,amssymb,latexsym,amsfonts,mathrsfs}
\usepackage{graphicx}
\usepackage{color}

\usepackage{fancyhdr}
\usepackage[top=30truemm, bottom=30truemm, left=25truemm, right=25truemm]{geometry}

\newcommand{\ep}{\varepsilon}
\newcommand{\nn}{\nonumber}

\newcommand{\SCR}[1]{{\mathscr #1}}

\newcommand{\CAL}[1]{{\cal #1}
}

\newcommand{\J}[1]{\left\langle #1 \right\rangle}
\newcommand{\D}[1]{{\mathscr D}( #1 )}

\theoremstyle{definition}
\newtheorem{theorem}{{\bf Theorem}}[section]

\newtheorem{lemma}[theorem]{{\bf Lemma}}
\newtheorem{proposition}[theorem]{{\bf Proposition}}

\newcounter{Exami}

\begin{document}

%%%%%%%%%%%%%%%%%%%%%%%%%%%%%%%%%%%%%%%%%%%%%%%%%%%%%%%%%%%%%%%%%%%%%%%%%%%%%%%%%%%%%%%%%%%%%%%%%%%%%%%%%%%%%%%%%%%%%%%%%%%
%author
%%%%%%%%%%%%%%%%%%%%%%%%%%%%%%%%%%%%%%%%%%%%%%%%%%%%%%%%%%%%%%%%%%%%%%%%%%%%%%%%%%%%%%%%%%%%%%%%%%%%%%%%%%%%%%%%%%%%%%%%%%%
\begin{flushleft}
{\bf \Large Strichartz estimates for quadratic repulsive potentials 
 } \\ \vspace{0.3cm} 
by {\bf \large Masaki Kawamoto
 $^{1}$} and {\bf \large Taisuke Yoneyama $^{2}$ } \\
$^{1}$ Graduate School of Science and Engineering, Ehime University, 3 Bunkyo-cho Matsuyama, Ehime 790-8577, Japan. \\ 
Email: kawamoto.masaki.zs@ehime-u.ac.jp \\ 
$^{2}$ Mathematics Section, Kitasato University
College of Liberal Arts and Sciences, 1-15-1, Kitazato, Minami-ku, Sagamihara-shi Kanagawa 22-0373, Japan. \\ 
 Email: t.yone26@kitasato-u.ac.jp
\end{flushleft}
%%%%%%%%%%%%%%%%%%%%%%%%%%%%%%%%%%%%%%%%%%%%%%%%%%%%%%%%%%%%%%%%%%%%%%%%%%%%%%%%%%%%%%%%%%%%%%%%%%%%%%%%%%%%%%%%%%%%%%%%%%%
%abst
%%%%%%%%%%%%%%%%%%%%%%%%%%%%%%%%%%%%%%%%%%%%%%%%%%%%%%%%%%%%%%%%%%%%%%%%%%%%%%%%%%%%%%%%%%%%%%%%%%%%%%%%%%%%%%%%%%%%%%%%%%%

\begin{abstract}
Quadratic repulsive potentials $- \tau ^2 |x| ^2$ accelerate the quantum
 particle, increasing the velocity of the particle exponentially in
 $t$; this phenomenon yields fast decaying dispersive
 estimates. In this study, we consider the Strichartz estimates associated with this phenomenon. 
 First, we consider the free repulsive Hamiltonian, and prove that the 
 Strichartz estimates hold for every admissible pair $(q,r)$, which satisfies $1/q +n/(2r) \geq n/4$ with $q$, $r \geq 2$. 
 Second, we consider the perturbed repulsive Hamiltonian with a slowly decaying potential, such that $|V(x)| \leq C(1+|x|)^{-\delta}$ 
 for some $\delta >0$, and prove that the Strichartz estimate holds with the same admissible pairs for repulsive-admissible pairs. 
\end{abstract}

\section{Introduction}

In this study, we consider the quadratic repulsive Hamiltonian $H_0$, which is defined by 
\begin{align*}
H_0= -\Delta  - \tau ^2x^2
=-\sum_{j=1}^n\frac{\partial^2}{\partial x_j^2}-\tau^2\sum_{j=1}^n x_j^2, 
\end{align*}
where $x = (x_1,x_2,...,x_n) \in {\bf R} ^n $, and $\tau \in{\bf R}\setminus\{0\}$ is a constant. 
Here, $H_0$ is known to be essentially self-adjoint on $C_0^\infty({\bf R} ^n)$,
and we denote the unique self-adjoint extension by the same symbol, i.e., $H_0$.
The following estimates play an important role in
certain issues associated with the repulsive Hamiltonian;
\begin{align} \label{1}
& \left\| 
e^{-i(t-s)H_0} \phi 
\right\|_{L^{\infty}({\bf R}^n)} \leq C  |\sinh (\tau (t-s))|^{-n/2} \left\|
 \phi \right\|_{L^1({\bf R}^n)}, 
 \end{align}
 and 
 \begin{align}
  \label{24}
& \left\| 
e^{-i (t-s) H_0} \phi 
\right\|_{L^2({\bf R}^n)} = \left\| \phi \right\|_{L^2({\bf R}^n)},
\end{align}
for $\phi  \in \SCR{S}({\bf R}^n) $ and $t,s \in {\bf R}$.
These estimates are deduced by the generalized Mehler formula 
\begin{align} \nn 
& (e^{-i\sigma H_0} \phi )(x) \\ &= \left( \frac{\tau}{2 \pi i \sinh 
 (2 \tau \sigma)}\right)^{n/2} \int_{{\bf R}^n} e^{i \tau ((x^2+y^2) \cosh
 (2 \tau \sigma)-2x \cdot y)/(2\sinh (2 \tau \sigma))} \phi(y) dy , \label{mdk1}
\end{align}
(see H\"{o}rmander \cite{Ho} 
and (2.10) of the Bony-Carles-H\"{a}fner-Michel theory \cite{BCHM}) and the unitarity of $e^{-i(t-s)H_0}$ on $L^2({\bf R}^n)$, respectively. 
The first theorem in this study involves the deduction of the Strichartz estimates associated with \eqref{1} and \eqref{24}; 
 
 Here, let us consider a case where $|t-s| \leq 1$. If $\kappa
 \geq n/2 $, then
\begin{align*}
|\sinh (\tau (t-s))^{-n/2} | \leq C |t-s|^{-n/2} \leq C_{\kappa} |t-s|^{-
 \kappa} ,
\end{align*}
for all $t,s \in {\bf R}$. 
Conversely, consider a case in which $|t-s| >1$. Then, for all $\kappa \geq
n/2 $,
\begin{align*}
|\sinh (\tau (t-s)) ^{-n/2}| \leq |t-s|^{- \kappa} |t-s|^{\kappa}
 |\sinh (\tau (t-s))|^{-n/2} \leq C |t-s|^{- \kappa}
\end{align*}
holds. Consequently, for all $\kappa \geq n/2$, we obtain
\begin{align}\label{2}
\left\|e^{-i(t-s) H_0} \phi (s) \right\|_{L^{\infty} ({\bf
 R}^n)} \leq C |t-s|^{- \kappa } \left\| \phi(s) \right\|_{L^1({\bf R}^n)}. 
\end{align} 
\begin{figure}[h]
\begin{minipage}{0.5\hsize}~~ \vspace{1.5cm}
  \begin{center}
   \includegraphics[width=5cm, bb= 0 0 600 200]{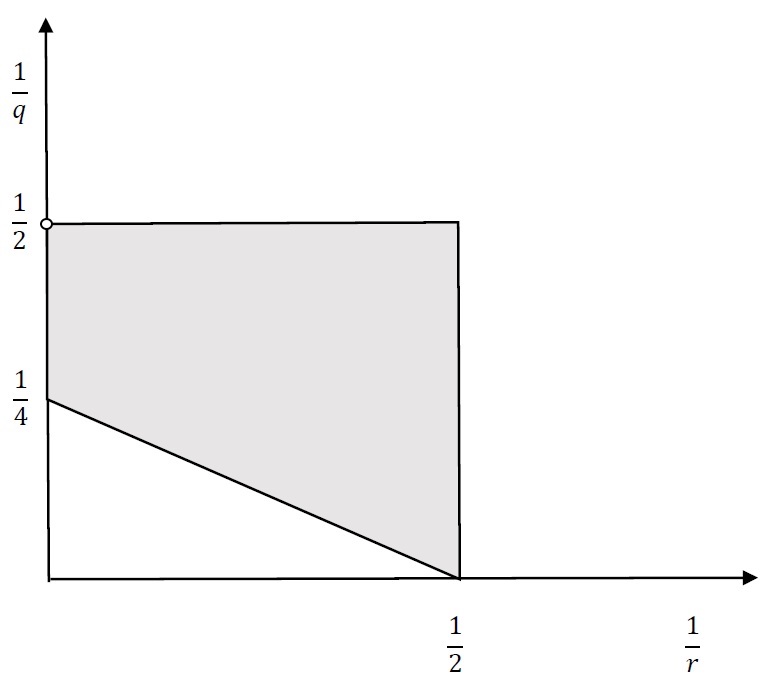}
  \end{center}
 \end{minipage}\qquad 
 \begin{minipage}{0.5\hsize}~~ \vspace{0.2cm}
  \begin{center}
  ~~\\ ~~ \\ ~~ \\ ~~ \\ ~~ \\ ~~ \\ 
   \includegraphics[width=5cm, bb= 0 0 600 200]{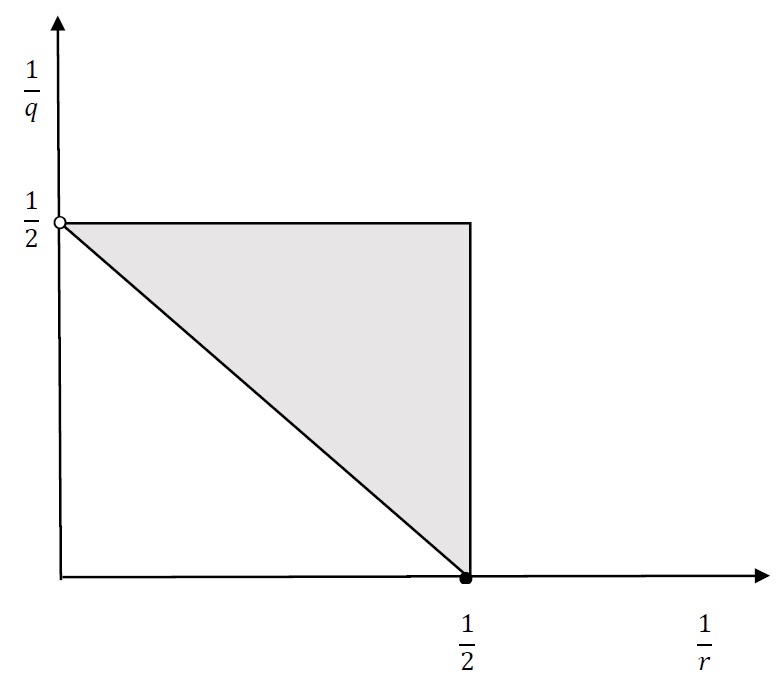}
  \end{center}
  \caption{  For $n=2$, the closed region comprises repulsive-admissible pairs}
  \label{f2}
 \end{minipage}
 \end{figure} 
 \begin{figure}[h]
 \begin{minipage}{0.5\hsize}
   \begin{center}
   \includegraphics[width=8cm, bb= 0 0 600 400]{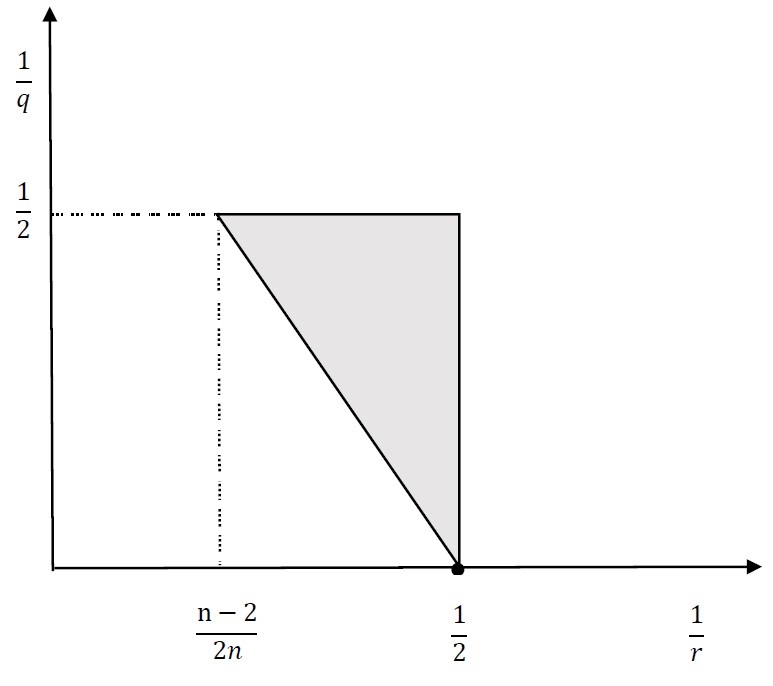}
  \end{center}
  \caption{For $n\geq3$, the closed region comprises repulsive-admissible pairs}
  \label{f3}
 \end{minipage}
 \end{figure}

\newpage

The dispersive estimates \eqref{2} and the consequence of Theorem 1.2 in Keel-Tao \cite{KT} yield the Strichartz estimates for all pairs $(q,r)$ satisfying $(q,r,\kappa) \neq (2, \infty, 1)$ and that 
\begin{align}\label{3}
\frac{1}{q} + \frac{\kappa }{r} = \frac{\kappa}{2}, \quad  q \geq 2 ,
\end{align}
for some $\kappa \geq n/2$. This proves the following Theorem \ref{T1}. Taking into account  
\begin{align*} 
\frac{1}{q} = \kappa  \left( \frac{1}{2} - \frac{1}{r} \right) \geq
 \frac{n}{2} \left(  \frac{1}{2} - \frac{1}{r} \right), 
\end{align*}
if pair $(q,r)$ satisfies 
\begin{align}\label{4}
\frac{1}{q} + \frac{n}{2r} \geq \frac{n}{4}, \quad  q \geq 2, \quad 2 \leq r \leq \frac{2n}{n-4/q},
\end{align}
with $(q,r,\kappa, n) \neq (2, \infty,1,2)$, $(2, \infty,1,1)$ then we call pair $(q,r)$ a {\em repulsive-admissible pair}.  For $n \geq 3$, pair $(q,r) = (2, 2n/(n-2))$ is referred to as the {\em end point} for the free Schr\"{o}dinger equation; 
the proof of the Strichartz estimates for the end point requires a more complicated argument. However, by virtue of the estimates, \eqref{1} and \eqref{24}, $(q,r) = (2, 2n/(n-2))$ can be included as the repulsive-admissible pair if $\tau >0$.

\begin{theorem}\label{T1}
Let $(q,r)$ and $(\tilde{q}, \tilde{r})$ be repulsive-admissible pairs.
Then, for all $f \in L^2({\bf R}^n)$ and $F \in L^{\tilde{q}'} ({\bf R} \, ; \, L^{\tilde{r}'} ({\bf R}^n))$, the following estimates hold.
\begin{align} \label{17}
\left\| 
e^{-itH_0} f 
\right\|_{L^q_tL^r_x} \leq C \left\| f \right\|_{L^2}
\end{align}
and 
\begin{align} \label{18}
\left\| 
\int_0^t e^{i(t-s)H_0} F(s) ds 
\right\|_{L^q_t L^r _x} \leq C \left\| F \right\|_{L^{\tilde{q}'}_t L^{\tilde{r}'}_x}
\end{align}
where 
\begin{align*}
\left\| F \right\|_{L_t^qL_x^r} := \left( \int\left\| F(s, \cdot )
 \right\|^q_{L^r({\bf R}^n)} ds \right)^{1/q}, \quad \left\| u
 \right\|_{L^r} = \left\|u \right\|_{L^r({\bf R}^n)}, 
\end{align*}
and $\tilde{q}'$ and $\tilde{r}'$ denote the H\"{o}lder conjugate of $\tilde{q}$ and $\tilde{r}$, i.e., $1/\tilde{q} + 1/\tilde{q}' =1$ and 
$1/\tilde{r} + 1/\tilde{r}' =1$.
\end{theorem}
  Thus, the Strichartz estimates hold in region $1/q+n/(2r) \geq n/4$ and $q,r \geq 2$ (see figures 1, 2, and 3); this type of estimate
has probably been unknown until now. Hence, exponential decay \eqref{1} expands the region
of $(q,r)$, on which the Strichartz estimates hold from line \eqref{3} with $\kappa = n/2$ to the region \eqref{4}. 

We apply Theorem \ref{T1} to the Hamiltonian $H=H_0 +V(x)$ with real valued $C^2$-function $V(x)$ and prove that the estimate \eqref{17} is satisfied for repulsive-admissible pairs by $e^{-itH}$ if $V(x)$ satisfies
\begin{align} \label{5}
|\nabla _x^{\alpha} V (x)| \leq C \J{x}^{-\delta },
\end{align} 
for $|\alpha| \leq 2$ with an arbitrarily small $\delta >0$, where $\J{\cdot} = (1 + |\cdot|^2)^{1/2}$.
Here, we can choose $\delta >0$ that is arbitrarily small. The second theorem in this paper is as follows.
\begin{theorem} \label{T2}
Let $(q,r)$ be a repulsive-admissible pair, $H=H_0 +V$ and $V(x)$ satisfies \eqref{5} with $\delta >0$.
Then, {\em for all} $f \in L^2({\bf R}^n)$, the following estimates hold:
\begin{align} \label{20}
\left\| e^{-it H} f \right\|_{L_t^{q} L^{r}_x}
 \leq C \left\| f \right\|_{L^2}.
\end{align}
\end{theorem}

In the case of $\tau =0$ and $V \equiv 0$, the Strichartz estimates were first considered and proven by Strichartz \cite{St} for the case of $q=r$. Further, the 
conditions of the admissible pair stated in \cite{St} were relaxed in some studies (e.g., Genibre-Velo \cite{GV} and Yajima \cite{Ya2}). In particular, 
the so-called {\em end-point Strichartz estimates} were proved \cite{KT}. Moreover, in case of $\tau = 0$ and $V \in L^{n/2} ({\bf R}^n )$, 
the Strichartz estimates were considered by Burq-Planchon-Stalker-Tahvildar Zadeh \cite{BPSZ2} (time-dependent case were studied by D'Ancona-Pierfelice-Visciglia \cite{DPV}, Goldberg \cite{G}, Rodonianski-Schlag \cite{RS} among others.) In here the condition  $V \in L^{n/2} ({\bf R}^n )$ is necessary condition for the time-in-global Strichartz estimate holds. Indeed, for $V(x)$ such that $|V(x)| \leq C\J{x}^{-2 - \delta _0}$, if $\delta _0\leq0$, Goldberg-Vega-Visciglia  \cite{GVV} proved that the Strichartz estimates did not always hold. After, Mizutani \cite{M2} relaxed the necessary condition of potentials for the Strichartz estimate to hold, in which he succeeded to include the potential so that $|\nabla^{\alpha} V(x)| \leq C\J{x}^{ - \delta _0 -|\alpha| }$ with $|\alpha| \leq n$ and $\delta _0 >0$, where $\delta _0$ can be taken arbitrary small, however he needed repulsive conditions $V(x) \geq C_1 (1+|x|)^{- \delta _0} $ and $- x \cdot \nabla V(x) \geq C_2 (1+|x|)^{-\delta_0} $, $C_1,C_2 >0$ in addition (or $V(x) = Z|x|^{- \gamma} + (\mbox{small perturbations})$ with $n \geq 3$, $Z>0$ and $\gamma \in (0,2)$). 

As the concrete model of a potential, which does not satisfy the condition of \cite{GVV}, \cite{M2} and is not in $L^{n/2} ({\bf R}^n)$,
the repulsive potential $V(x) = C_{\mathrm{ha}} |x|^{-2}$ has been studied by several researchers. 
They obtained various results related to the Strichartz estimates
(e.g., Burq-Planchon-Stalker-Tahvildar Zadeh \cite{BPSZ} \cite{BPSZ2}, Mizutani \cite{M},
and Pierfelice \cite{P}), where $C_{\mathrm{ha}}$ is a positive constant that is smaller than or equal to the best constant $-(n-2)^2/4$ of Hardy's.

Recently, certain issues associated with the quadratic repulsive Hamiltonian have been studied by Bony-Carles-H\"{a}fner-Michel
\cite{BCHM}, Carles \cite{C}, Fang-Han \cite{FH}, Ishida \cite{Is} among others. In particular, Carles \cite{C} studied certain problems related to 
the nonlinear equation, $i \partial _t v = H_0 v + \lambda|v|^{2b} v$, 
with $b < 2/(n-2)$ and $\lambda \in {\bf R}$ using the Strichartz estimates with the admissible pair $(q,r)$, which satisfies $1/q + n/ (2r) = n/4$ which is a part of \eqref{4}. 
Hence, our results may allow us to improve some results of \cite{C} and enable us to consider the nonlinear Schr\"{o}dinger equations with quadratic repulsive potential $-\tau ^2 x^2$ and perturbation $V$.

Now, let us discuss the case of $\tau \in {\bf C}$. It is well known that, if $\tau ^2 $ is negative, i.e., $\tau \in{\bf C}$ with $\mathrm{Re }\tau = 0$,
it is impossible to prove the global-in-time Strichartz estimates.
This is because term $|\sinh (\tau (t-s))|^{-n/2} $ in the dispersive estimates never converges to $0$ as $|t-s| \to \infty$.
Hence, only the local-in-time Strichartz estimates (e.g., Yajima-Zhang \cite{YZ}) were proven. 

Conversely, if $\tau ^2$ is positive, coefficient $|\sinh (\tau (t-s))|^{-n/2} $ converges to $0$ 
exponentially, enabling us to prove the global-in-time Strichartz estimates.
Here, we remark that the dispersive estimate \eqref{1} is ``inhomogeneous,'' that is,
the coefficient of \eqref{1} acts as $|\sinh (\tau (t-s))|^{-n/2}$ for $|t-s| \geq 1$; however, it acts as $|t-s|^{-n/2}$ for $|t-s| \leq 1$. 
In this case, it may be difficult to apply the method of Keel-Tao \cite{KT} directly.
Kawamoto-Yoneyama \cite{KY} considered the Strichartz estimates
for inhomogeneous dispersive estimates; more precisely, they considered a case wherein the coefficient of dispersive estimate acts as $|t-s|^{-n(1-\varepsilon)/2}$ 
for $|t-s| \geq 1$, as well as $|t-s|^{-n/2}$ for $|t-s| \leq 1$ and for some $0< \varepsilon <1/2$. 
In this case, by introducing the time-weighted space, we succeeded in avoiding the difficulties arising from inhomogeneous dispersive estimates. 
In the quadratic repulsive case, to overcome this difficulty, 
 we reduced the decaying order from 
$|\sinh (\tau (t-s))|^{-n/2}$ to $|t-s| ^{-\kappa}$ for $|t-s| \geq 1$ with $\kappa \geq n/2$.
It is important to use \eqref{2} instead of \eqref{1} to apply Keel-Tao's result and obtain Theorem \ref{T1} in spite of that \eqref{2} is much weaker than \eqref{1} and, that they were unable to deal with potentials which decays like $(\log(1+|x|))^{-2 - \ep}$ because of this, where it may be expected that we can prove the Strichartz estimate under such log-decay potential if the quadratic repulsive potential exists.

\section{Resolvent estimates and proof of theorem \ref{T2}}

Before showing Theorem \ref{T2}, we introduce the scheme for the proof. Here, we mimic the approach of \cite{RS} (see also \cite{BPSZ}).
By Duhamel's formula, we obtain 
\begin{align*}
e^{-itH} u = e^{-itH_0}u -i \int _0^t e^{-i(t-s) H_0} V e^{-isH} u ds, 
\end{align*}
for $u \in L^2({\bf R}^n)$.
We let $\mu = \mu_{\delta} > \max\{ 2,n, 2n/ \delta \} $. Then we have $|V| ^{1/2} \in L^{\mu} ({\bf R}^n)$ and obtain 
\begin{align*}
\left\| V e^{-itH} u\right\|_{L^{2 \mu /(\mu + 2)}} \leq C \left\| |V|
 ^{1/2}\right\|_{L^{\mu}} \left\| |V|^{1/2} e^{-itH} u \right\|_{L^2}. 
\end{align*} 
Because $\mu > 2$  and $\mu > n$, the H\"{o}lder conjugate of $(2 \mu/(\mu +2))$ satisfies 
$(2 \mu /(\mu
+2)) ' = 2 \mu /(\mu -2) > 0$ and
\begin{align*}
\frac12 + \frac{n(\mu -2)}{4 \mu} =\frac12 - \frac{n}{2 \mu} +
 \frac{n}{4} \geq \frac{n}{4},
\end{align*}
which implies that $(q,r) = (2, 2 \mu/(\mu -2))$ is a repulsive-admissible pair. By using Theorem \ref{T1}, we have 
\begin{align} \nn
\left\| 
e^{-itH}u
\right\|_{L_t^{q} L^{r}_x} & \leq C \left\| u \right\|_{L^2} + C \left\| V
 e^{-itH}u \right\| _{L^2_t L^{2 \mu/(\mu +2)}_x} \\ 
& \leq C \left\| u \right\|_{L^2} + C \left( 
\int_{\bf R} \left\| |V|^{1/2} e^{-itH} u \right\|_{L^2} ^2 dt 
\right)^{1/2}.
\label{ad2}
\end{align} 
Hence, if one obtains 
\begin{align} \label{7}
 \int_{\bf R} \left\| |V|^{1/2} e^{-itH} u \right\|_{L^2} ^2 dt 
\leq C \left\| u \right\|_{L^2}^2,
\end{align}
the proof of Theorem \ref{T2} is completed.

In the following, we prove \eqref{7}.
Here, we use the Kato smooth perturbation theory (see Kato \cite{Ka}).
First,
the following lemma holds.
 \begin{lemma}[D'Ancona \cite{D}] \label{L3}
 Let $G$ be a closed operator with dense domain $\D{G} \subset L^2({\bf R}^n)$. Assume that there exist $C_K >0$ and $\nu_0 >0$, such that 
 \begin{align} \label{27}
\gamma _1:= \sup_{\lambda \in {\bf R}, \, |\nu| \leq \nu_0, }\left| \left( G (H- \lambda -i \nu)^{-1} G^{\ast}u , u \right)_{L^2({\bf R}^n)}
 \right| \leq C_{\mathrm{K}} \left\| u \right\|_{L^2}^2
 \end{align} 
holds for any $u \in L^2({\bf R}^n)$. Then, the following inequalities hold:
\begin{align*} 
\gamma _2 &:=  \int_{\bf R} \left\| G e^{-itH} u \right\|^2_{L^2}
 dt \leq  4C_{\mathrm{K}} \left\| u \right\|_{L^2}^2 
 \end{align*}
  and 
  \begin{align*}
\gamma _3 &:=  \int_{\bf R} \left\| e^{- |\nu| t} \int_0^t Ge^{-i(t-s)H} G^{\ast} h(s) ds \right\|_{L^2} ^2 dt \leq 4C_{\mathrm{K}}^2  
\int_{\bf R} \left\|  e^{- |\nu| t} h(s) \right\|_{L^2} ^2 ds 
\end{align*} 
for $
\nu \in [-\nu_0, \nu_0] $.
\end{lemma}
Here, we replace $G$ with $|V|^{1/2}$ in \eqref{27}, and prove $\gamma _1 \leq C_{\mathrm{K}} \left\| u \right\|_{L^2}^2$. 
Then, using Lemma \ref{L3} with $ G=|V|^{1/2}$, we have \eqref{7}, which completes the proof of Theorem \ref{T2}.
Here, it is easily demonstrated that 
\begin{align*}
\gamma_1 \leq \sup_{\lambda \in {\bf R}, \, |\nu| \leq \nu_0 } \left\| G (H- \lambda -i \nu)^{-1} G^{\ast} \right\|_{\SCR{B}(L^2({\bf R}^n))} \left\| u \right\|_{L^2}^2 
\end{align*}
holds, where $\left\| \cdot \right\|_{\SCR{B}(L^2({\bf R}^n))}$ is the operator norm on $L^2({\bf R}^n)$. Hence, in the sequel, we prove the so-called uniform resolvent estimate 
\begin{align} \label{26}
\sup_{\lambda \in {\bf R}, \, 0< \nu \leq \nu _0}\left\|  |V|^{1/2} (H- \lambda  \mp i \nu)^{-1} |V|^{1/2} \right\|_{\SCR{B}(L^2({\bf R}^n))} \leq C,
\end{align}
for some small positive constant $\nu_0 >0$.

The proof is divided into two parts. In the first part, we prove \eqref{26} for $\lambda$ in compact interval $[-R,R]$. The key estimate is the weighted resolvent estimate of Theorem 4.1 of \cite{BCHM}. Here, we note that, if $V$ satisfies \eqref{5}, 
it immediately holds $V \in L^p ({\bf R}^n)$ for large $p >1$.
Hence, by virtue of Theorem 2.8 and Remark 2.9 of \cite{BCHM} or Corollary 1.3 of Itakura \cite{It}, $H$ has no eigenvalues. Then, Theorem 4.1 of \cite{BCHM} proves the limiting absorption principle; for some $s > 1/2$ and a small interval $\Lambda \subset {\bf R}$, there exists $C_{\Lambda} > 0$ such that
 \begin{align}\label{ad1}
 \sup_{\lambda \in \Lambda, \mu >0 } \left\| \J{\SCR{A}}^{-s} (H - \lambda \mp i \mu)^{-1} \J{\SCR{A}}^{-s} f \right\|_{L^2} \leq C_{\Lambda} \| f \|_{L^2} 
 \end{align}
 with a self-adjoint operator 
  \begin{align*}
 \SCR{A} = 2 |\tau| \log(1 + (-i \nabla + x)^2 ) - 2|\tau| \log(1 + (-i \nabla - x)^2),
 \end{align*}
 which is defined as a {\em pseudo-differential operator}, see (3.2) of \cite{BCHM}. Hence, the aim of the first part is to extend $\Lambda$, a small interval, to a large interval $[-R,R]$ (Proposition \ref{T3}), and replace the weight in \eqref{ad1} from $\J{\SCR{A}}^{-s}$ to $\J{x}^{-\rho}$ with some $\rho >0$ (Proposition \ref{T4}). 
 
In the second part, we provide the proof of \eqref{26} for the high-energy case using the argument in Yajima \cite{Ya} and the resolvent estimate as
\begin{align*} 
\lim_{\lambda \to \infty} \sup_{0 <\nu \leq \nu_0} \left\| |V|^{1/2} (H_0- \lambda
 \mp i \nu )^{-1} |V^{1/2} | \phi \right\|_{L^2} =0.
\end{align*}
%%%%%%%%%%%%%%%%%%%%%%%%%%%%%%%%%%%%%%%%%%%%%%%%%%%%%%%%%%%%%%%%%%%%%%%%%%%%%%%%%%%%%%%%%%%%%%%%%%%%%%%%%%%%%%%%%%%%%%%%%%%%%%%%%%%%%%%%%%%%%%%%%%%%%%%%5

\subsection{Low-energy estimate}

First, we proved the low-energy resolvent estimate. Here, we prove the following proposition. 
\begin{proposition}\label{T3}
Let $V \in C^2({\bf R}^n)$ satisfy \eqref{5}, $R>0$ and $s>1/2$. Then there exists constant $C_{R,1}>0$ such that for all $f \in L^2({\bf
 R}^n)$, 
\begin{align}\label{8}
  \sup_{\lambda \in [-R,R], \, 0< \nu \leq \nu_0
 } \left\| \J{\SCR{A}} ^{-s} (H -
 \lambda  \mp i \nu )^{-1} \J{\SCR{A}}^{-s} f \right\|_{L^2} \leq C_{R,1} \left\| f
 \right\|_{L^2}.
\end{align} 
\end{proposition} 
\begin{proof} 
By Theorem 4.1 of \cite{BCHM} and the fact that $\sigma _{\mathrm{pp}} (H)= \emptyset $, there exists sequence $\left\{ \theta _{j} \right\}_{j \in \{0,1,\ldots, J-1\}}$, $J \in {\bf N}$, and $J< \infty$ with
$\theta _j - \theta _{j-1} >0$ such that
\begin{align*}
\sup_{\lambda \in  \Gamma_j, \, 0< \nu \leq \nu_0
 } \left\| \J{\SCR{A}} ^{-s} (H -
 \lambda  \mp i \nu )^{-1} \J{\SCR{A}}^{-s} f \right\|_{L^2} \leq C_{j} \left\| f
 \right\|_{L^2},
\end{align*}
where $\Gamma _j = [\theta _{j-1} , \theta _j]$ and $C_j>0$ is a constant, and depends only on $\nu_0$ and $\theta _j - \theta _{j-1}$.
Hence, by allowing $[-R, R] \subset \bigcup_{j=1}^{j=J} \Gamma _j$, we have 
\begin{align*}
& \sup_{ \lambda \in [-R,R], \, 0< \nu \leq \nu_0} \left\| \J{\SCR{A}} ^{-s} (H -
 \lambda  \mp i \nu )^{-1} \J{\SCR{A}}^{-s} f \right\|_{L^2} \\ & \leq \sum_{j =1}^J \sup_{ \lambda \in [\theta_{j-1},\theta_j], \, 0< \nu \leq \nu_0}  
 \left\| \J{\SCR{A}} ^{-s} (H -
 \lambda  \mp i \nu )^{-1} \J{\SCR{A}}^{-s} f \right\|_{L^2} \\
 & \leq \sum_{j =1}^J C_j \left\| f \right\|_{L^2} \leq C_{R,1} \left\| f \right\|_{L^2} .
\end{align*}
\end{proof}

The following theorem can be obtained as a consequence of Proposition \ref{T3}.
\begin{proposition}[Low-energy estimate] \label{T4}
Under the same assumptions of Proposition \ref{T3}, for all $f \in L^2({\bf R}^n)$, 
\begin{align} \label{21}
 \sup_{\lambda \in [-R,R], \, 0< \nu \leq \nu_0
 } \left\| \J{x}^{-\delta /2} (H -
 \lambda  \mp i \nu )^{-1} \J{x}^{-\delta /2} f \right\|_{L^2} \leq C_{R,2} \left\| f
 \right\|_{L^2}
\end{align}
holds, where $C_{R,2}$ depends on $R$.
\end{proposition}
\begin{proof} 
For fixed $\lambda \in [-R,R]$ and $\mu \in (0, \mu_0]$, we have
\begin{align*}
& \left\| \J{x}^{-\delta /2} (H -
 \lambda  \mp i \nu )^{-1} \J{x}^{-\delta /2} f \right\|_{L^2} \\ & \leq \left\| \J{x}^{- \delta /2} \J{\SCR{A}}^{s} \right\|_{\SCR{B}(L^2({\bf R}^n))} ^2 \left\| 
 \J{\SCR{A}}^{-s} (H - \lambda  \mp i \nu )^{-1} \J{\SCR{A}}^{-s} \right\|_{\SCR{B}(L^2({\bf R}^n))} \left\| f \right\|_{L^2}.
\end{align*}
Here, according to (3.13), (3.14), and (3.15) of \cite{BCHM}, $\J{x}^{-\delta /2} \J{\SCR{A}}^{s}$ is a bounded operator.
Thus, together with \eqref{8}, inequality \eqref{21} is obtained.
\end{proof}

Unfortunately, because the right-hand side of \eqref{21} depends on $R$ and the limiting absorption principle fails in the case of $R \to \infty$, we need to provide
another estimate for the high-energy case. 

\subsection{High-energy estimates}
We define 
\begin{align*}
R_{\rho} (\pm \nu) = \J{x}^{- \rho} (H_0- \lambda \mp i \nu)^{-1}
 \J{x}^{- \rho }, 
\end{align*} 
where $\rho >0$ is a constant that is given later. For a fixed $\phi \in L^2({\bf R}^n)$, we define
\begin{align*} 
\CAL{A}_0 &:= \left\| R_{\rho} (\pm \nu) \phi \right\|_{L^2} = \left\|  \J{x}^{-0} R_{\rho}(\pm \nu)
 \J{x}^{-0} \phi \right\|_{L^2}, \\ 
\CAL{A} _2 &:= \left\| \J{x}^{-2} R_{\rho}(\pm \nu) \J{x}^{-2} \phi \right\|_{L^2}.
\end{align*}
In this subsection, we prove that, for some constants $C_0>0$ and $C_2 >0$,
which are independent of $\nu$, $\CAL{A}_0 \leq C_0 \left\| \phi \right\|_{L^2}$ and $\CAL{A}_2 \leq C_2 
|\lambda|^{-1/3} \left\| \phi \right\|_{L^2}$ hold. Then, based on the interpolation
theorem that has arisen from the Hadamard three lines theorem, it
follows that for $0 \leq \theta \leq 2$, 
\begin{align} \label{13}
\left\| 
\J{x}^{- \theta} R_{\rho} (\pm \nu) \J{x}^{- \theta}\phi
\right\|_{L^2} \leq C_0^{1-\theta/2} C_2^{\theta/2} |\lambda|^{-\theta /6}
 \left\| \phi \right\|_{L^2}
\end{align}
holds. By using this estimate, we can prove \eqref{26} for the high-energy
case. To deduce the estimates of $\CAL{A}_0$ and $\CAL{A}_2$, we
introduce the following lemma. 
%%%%%%%%%%%%%%%%%%%%%%%%%%%%%%%%%%%%%%%%%%%%%
\begin{lemma} \label{L1}
Let $\rho >0$ and $Q \geq 2$ with $\rho Q 
 >n$ and $Q > n$. 
Then, there exists constant $C>0$, such that 
\begin{align} \label{11}
\int_0^{\infty} \left\| 
\J{x}^{- \rho} e^{-i \sigma H_0} \J{x}^{- \rho} \phi 
\right\|_{L^2}  d \sigma \leq C\left\|\J{\cdot}^{- \rho} \right\|^2_{L^Q} \left\| \phi \right\|_{L^2}.
\end{align}
\end{lemma}
\begin{proof}
We mimic the approach of Kato \cite{Ka}. By using the H\"{o}lder inequality for $1/P + 1/Q = 1/2 $, we have 
\begin{align} \label{X1}
\left\| 
\J{x}^{- \rho} e^{-i \sigma H_0} \J{x}^{- \rho} \phi 
\right\|_{L^2} \leq C  \left\| \J{\cdot}^{- \rho} \right\|_{L^{Q}} 
\left\| 
e^{-i \sigma H_0} \J{x}^{- \rho} \phi 
\right\|_{L^{P}}.  
\end{align}
By using \eqref{1}, \eqref{24} and Riesz-Tholin's interpolation theorem 
for $P \geq 2$ and $v \in L^{P/(P-1)}({\bf R}^n)$, we have 
\begin{align} \label{X2}
\left\| 
e^{-i \sigma H_0} v
\right\|_{L^P} \leq C |\sinh (2\tau \sigma)|^{-n(1/2-1/P)}\left\| v
 \right\|_{L^{P/(P-1)}}.
\end{align}
Hence, \eqref{X1} and \eqref{X2} yield
\begin{align}\label{12}
\left\| 
\J{x}^{- \rho} e^{-i \sigma H_0} \J{x}^{- \rho} \phi 
\right\|_{L^2} \leq C  |\sinh (2\tau \sigma)|^{-n/Q} \left\|\J{\cdot}^{-
 \rho} \right\|^2_{L^Q}  \left\| \phi \right\|_{L^2}.  
\end{align}
This inequality yields 
\begin{align*}
(\mbox{l.h.s. of } \eqref{11}) \leq C \left\| \J{\cdot}^{-
 \rho}\right\|_{L^Q}^2 \left\| \phi \right\|_{L^2}\int_0^{\infty} |\sinh (2 \tau \sigma)|^{-n/Q} d
 \sigma ,
\end{align*}
which yields the proof of Lemma \ref{L1}.
\end{proof}
Now, we prove $\CAL{A}_2 \leq C_2 |\lambda|^{-1/3} \left\| \phi
\right\|_{L^2}$. Inequality $\CAL{A}_0 \leq C_0$ holds by virtue of Lemma \ref{L1} and the Laplace transform. 

We only consider the case of $\lambda >0$ and
$0< \nu \leq  \nu_0$ because the other cases can be proven in a similar way. 
By using the Laplace transform, 
\begin{align*}
\CAL{A}_2 = \left\| 
\int_0^{\infty} \J{x}^{- \rho-2} e^{-i \sigma (H_0 - \lambda -i \nu)}
 \J{x}^{- \rho -2} \phi d \sigma
\right\|_{L^2} .
\end{align*}
We split the interval by splitting the integral $(0, \infty) = (0, \lambda ^{-1/3}) \cup [\lambda^{-1/3} , \infty ) $. For the integral over $(0, \lambda^{-1/3} )$, we apply $\CAL{A}_0 \leq C$ and estimate it by 
\begin{align*}
I_1 = \left\| \int_0^{\lambda ^{-1/3}} \J{x}^{- \rho-2} e^{-i \sigma (H_0 - \lambda -i \nu)}
 \J{x}^{- \rho -2} \phi d \sigma\right\|_{L^2} \leq  C \left\| \phi\right\|_{L^2} \int_0^{\lambda ^{-1/3}} d \sigma \leq
 C \lambda ^{-1/3} \left\| \phi \right\|_{L^2}.
\end{align*}
For integral over $[\lambda ^{-1/3}, \infty ) $, we multiply the integrand by $1 = - \lambda^{-1}((H_0 - \lambda -i \nu) - (H_0 -i \nu)) $, integrate out the term $(H_0- \lambda -i \nu)e^{-i \sigma (H_0- \lambda -i \nu)} = i (d/d \sigma) e^{-i \sigma (H_0 - \lambda - i\nu)} $. Then using \eqref{11}, we estimate it by
\begin{align*}
I_2 &= \frac{1}{\lambda} \left\| 
\J{x}^{-2- \rho}\left[ 
e^{-i \sigma (H_0 - \lambda -i \nu)}
\right]^{\sigma = \infty}_{\sigma = \lambda^{-1/3}} \J{x}^{- 2 -
 \rho} \phi 
\right\|_{L^2} \\ & \quad + 
\frac{1}{\lambda} \left\| \int_{\lambda^{-1/3}}^{\infty} \J{x}^{- \rho -2} (H_0 -i\nu)e^{-i \sigma (H_0 - \lambda -i \nu)}
 \J{x}^{- \rho -2} \phi d \sigma
\right\|_{L^2} \\ 
& \leq C \nu_0 |\lambda|^{-1
} \left\| \phi \right\|_{L^2} + \frac{1}{\lambda} \left\| \int_{\lambda^{-1/3}}^{\infty} \J{x}^{-2- \rho} H_0 e^{-i \sigma (H_0 - \lambda -i \nu)}
 \J{x}^{- \rho -2} \phi d \sigma
\right\|_{L^2}.
\end{align*}
This inequality follows from \eqref{12}. Since $H_0 e^{-i \sigma H_0} = i (\partial / \partial \sigma)e^{-i \sigma H_0} $, differentiating the integral kernel \eqref{mdk1} of $e^{-i \sigma H_0}$ by $\sigma$, we obtain
\begin{align*}
 (\J{x}^{-2} H_0 e^{-i\sigma H_0} \J{x}^{-2} \phi)(x) = A_1 + A_2 + A_3 +A_4, 
\end{align*} 
where
\begin{align*} 
A_1 &=( \J{x}^{-2} x^2 (\tau/ \tanh(2\tau \sigma))^2 e^{-i \sigma H_0}
 \J{x}^{-2} \phi)(x), \\ 
A_2 &= n( \J{x}^{-2} (-i\tau / \tanh(2\tau \sigma)) e^{-i \sigma H_0}
 \J{x}^{-2} \phi)(x), \\ 
A_3 &= ( \J{x}^{-2}  (\tau / \sinh(2\tau \sigma))^2 e^{-i \sigma H_0}
 x^2 \J{x}^{-2} \phi)(x), \\
A_4 &= - \sum_{j=1}^n ( \J{x}^{-2} x_j (\tau^2  \cosh (2 \tau \sigma)/ \sinh(2\tau \sigma)^2) e^{-i \sigma H_0}
 x_j \J{x}^{-2} \phi)(x). 
\end{align*}
For $\sigma \in [\lambda ^{-1/3}, \infty)$, $|\sinh (2 \tau \sigma)| \geq C, 
\lambda ^{-1/3}$ holds, which gives 
\begin{align*}
& \left\|\int_{\lambda ^{- 1/3}}^{\infty} \J{x}^{-2- \rho} H_0 e^{-i
 \sigma (H_0 - \lambda -i \nu)} \J{x}^{-2 - \sigma} \phi d \sigma
 \right\|_{L^2} \\
& \leq C\lambda^{2/3} \left\| x^2 \J{x}^{-2} \right\|_{\SCR{B}(L^2({\bf R}^n))} 
\\ & \qquad \times\int_{0}^{\infty} \left\| \J{x}^{- \rho} e^{-i
\sigma H_0} \J{x}^{-
 \rho} \cdot ((x^2 +1)\J{x}^{-2} \phi) \right\|_{L^2} d \sigma \\ 
& \leq C \lambda^{2/3}\left\| \phi \right\|_{L^2}.
\end{align*}
Thus, we have $\CAL{A}_2 \leq I_1 + I_2 \leq C_2 \lambda
^{-1/3}$.
By using \eqref{13} with $\rho =
\theta = \delta /4$, we obtain
\begin{align} \label{14}
\sup_{ 0< \nu \leq \nu_0 } 
\left\| 
|V|^{1/2} (H_0- \lambda \mp i \nu)^{-1} |V|^{1/2} \phi 
\right\|_{L^2} \leq C_0^{1-\delta /8} C_2^{\delta /8} |\lambda|^{-
 \delta /24} \left\| \phi \right\|_{L^2}.
\end{align}
%%%%%%%%%%%%%%%%%%%%%%%%%%%%%%%%%%%%%%%%%%%%%%%%%%%%%%%
\begin{theorem}\label{T7}
Let us assume that $V\in C^2({\bf R}^n)$ satisfies \eqref{5}. Then, for all $\phi \in L^2({\bf R}^n)$, there exists $C>0$ such that
\begin{align*}
\sup_{|\lambda | \geq R, \, 0< \nu \leq \nu_0} \left\| 
|V|^{1/2} (H- \lambda \mp i \nu )^{-1} |V|^{1/2} \phi 
\right\|_{L^2} \leq C \left\| \phi \right\|_{L^2}.
\end{align*}
\end{theorem}
\begin{proof}
First, by considering \eqref{14}, we notice that, for large $\lambda \gg 1$, the Birman-Schwinger operator 
\begin{align*}
1 + |V|^{1/2} \mathrm{sign}(V) (H_0 - \lambda \mp i \nu)^{-1} |V|^{1/2}, 
\end{align*}
is invertible uniform in $\nu$ and $\lambda \gg 1$.
Then, by mimicking the approach of Yajima \cite{Ya}, it follows that 
\begin{align*}
& \rho _1 (H - \lambda \mp i \nu)^{-1} \rho_2  -
 \rho_1 (H_0 - \lambda \mp i \nu)^{-1} \rho_2 \\
 & = 
\rho_1 (H_0 - \lambda \mp i \nu)^{-1} \rho_2
\left(1 + \rho_1 (H_0 - \lambda \mp i \nu)^{-1}
 \rho_2 \right)^{-1}
\rho_1 (H_0 - \lambda \mp i \nu)^{-1} \rho_2
\end{align*}
with $\rho_1=  |V|^{1/2} \mathrm{sign}(V) $ and $ \rho_2 =|V|^{1/2}$,
which completes the proof.
\end{proof}
\subsection{Proof of theorem \ref{T2}}
Finally, we prove the main theorem. Combining Proposition \ref{T4} and Theorem \ref{T7}, one has the limiting absorption principle, as shown 
\begin{align*}
\sup_{\lambda \in {\bf R},\, 0 < \nu \leq \nu_0} \left\| |V|^{1/2} (H - \lambda \mp i \nu)^{-1} |V|^{1/2} \phi \right\|_{L^2} \leq C \| \phi \|_{L^2}.
\end{align*} 
This inequality clearly enables the following to be deduced, 
\begin{align*}
\sup_{\lambda \in {\bf R},\, |\nu | \leq \nu_0} \left| \left( 
|V|^{1/2} (H - \lambda - i \nu)^{-1} |V|^{1/2} \phi , \phi \right)_{L^2({\bf R}^n)}
\right| \leq C \| \phi \|_{L^2}^2. 
\end{align*}
Then, owing to the result of \cite{D} with $G = |V|^{1/2} $, we have 
\begin{align*}
\int_{\bf R} \left\||V|^{1/2} e^{-itH} \phi \right\|_{L^2}^2 \leq 4 C \| \phi \|^2_{L^2}. 
\end{align*} 
By considering \eqref{ad2}, Theorem \ref{T2} can be obtained.

% Non-BibTeX users please use

\end{document}